\documentclass[a4paper,12pt]{article}
\usepackage{amssymb}
\usepackage{amsthm}
\usepackage{amsxtra}
\usepackage{amsfonts}
\usepackage{xr}
\usepackage{color}
\usepackage{colordvi}
\usepackage{enumitem}
\usepackage{pb-diagram}%, pb-xy, xy}

\begin{document}
%\def\red{\textcolor{red}}
%\definecolor{green}{rgb}{0.1,.5,.1}
%\def\green{\textcolor{green}}
%\def\blue{\textcolor{blue}}
%\def\magenta{\textcolor{magenta}}
%\definecolor{grey}{rgb}{0.5,0.5,0.5}
%\def\grey{\textcolor{grey}}

\newcommand\red[1]{{\color{red}#1}}
\newcommand\blue[1]{{\color{blue}#1}}

\newcommand\sect{\section}

\newtheorem{thm}{Theorem}[section]
\newtheorem{cor}[thm]{Corollary}
\newtheorem{lem}[thm]{Lemma}
\newtheorem{prop}[thm]{Proposition}
\newtheorem{propconstr}[thm]{Proposition-Construction}
\newtheorem{tool}[thm]{Main tool/Sublemma}

\theoremstyle{definition}
\newtheorem{para}[thm]{}
\newtheorem{ax}[thm]{Axiom}
\newtheorem{conj}[thm]{Conjecture}
\newtheorem{defn}[thm]{Definition}
\newtheorem{notation}[thm]{Notation}
\newtheorem{rems}[thm]{Remarks}
\newtheorem{rem}[thm]{Remark}
\newtheorem{question}[thm]{Question}
\newtheorem{example}[thm]{Example}
\newtheorem{problem}[thm]{Problem}
\newtheorem{exercise}[thm]{Exercise}
\newtheorem{ex}[thm]{Exercise}
\newtheorem{fact}[thm]{Facts}

\overfullrule=0pt

\newcommand{\si}{\sigma}
\newcommand{\prf}{\smallskip\noindent{\it        Proof}. }
\newcommand{\call}{{\mathcal L}}
\newcommand{\nat}{{\mathbb  N}}

\newcommand{\inv}{^{-1}}

\newcommand{\trdeg}{\mathrm{tr.deg}}

\newcommand{\rest}{{\lower       .25     em      \hbox{$\vert$}}}

\newcommand{\zee}{{\mathbb  Z}}
\newcommand{\conc}{^\frown}
\newcommand{\acl}{\mathrm{acl}}

\newcommand{\cals}{{\mathcal S}}

\newcommand{\aut}{\mathrm{Aut}}
\newcommand{\ffi}{{\mathbb  F}}
\newcommand{\ffiti}{\tilde{\mathBbb          F}}

\newcommand{\gal}{\mathrm{Gal}}

\newcommand{\rat}{{\mathbb Q}}

\newcommand{\fix}{\mathrm{Fix}}
\newcommand\tp{\mathrm{tp}}

\newcommand{\calt}{{\mathcal T}}

\newcommand\vlabel[1]{\label{#1} {\rm[#1]}}

\title{A note on the non-existence of prime models of theories of
  pseudo-finite fields}

\author{Zo\'e Chatzidakis\thanks{partially supported  by
    ANR--13-BS01-0006 (ValCoMo) and ANR GeoMod (ANR-DFG, AAPG2019). Most of this work was done while the
    author was at the Institut Henri Poincar\'e, during the trimester
    {\em Model theory, combinatorics and valued fields}}\\ %% DMA (UMR 8553), Ecole Normale Sup\'erieure\\
  %% CNRS, PSL Research University}
  CNRS (IMJ-PRG), Sorbonne Universit\'e, Universit\'e Paris Cit\'e}

%\date{}
%\centerline{\today}
\maketitle

\begin{abstract}We show that if a field $A$ is not pseudo-finite, then there is no prime model of the theory of pseudo-finite fields over $A$. Assuming GCH, we extend this result to $\kappa$-prime models, for $\kappa$ an  uncountable cardinal or $\aleph_\varepsilon$.
  \end{abstract}
\section*{Introduction}
In this short note, we show that prime models of the theory of
pseudo-finite fields do not exist. More precisely, we consider the
following theory $T(A)$: $\ffi$ is a pseudo-finite field, $A$ a
relatively algebraically closed subfield of $\ffi$, and $T(A)$ is the
theory of the field $\ffi$ in the language of rings augmented by
constant symbols for the elements of $A$. Our first result is: \\[0.05in]
{\bf Theorem \ref{thm1}}. {\em Let $T(A)$ be as above. If $A$ is not
  pseudo-finite, then $T(A)$ has no prime model.}\\[0.05in]
When $A$ is infinite, the proof is done by constructing $2^{|A|}$ non-isomorphic
models of $T(A)$, of transcendence degree $1$ over $A$ (Proposition
\ref{prop1} and Remark \ref{rem-prop1}). \\[0.05in]
Next we address the question of existence of $\kappa$-prime models of
$T(A)$, where $\kappa$ is an uncountable cardinal or
$\aleph_\varepsilon$. We assume GCH, and again show in Theorem
\ref{thm2}  the non-existence of
$\kappa$-prime models of $T(A)$ (when $A$ is not already
$\kappa$-saturated pseudo-finite) in the following cases: if $\kappa\geq
\aleph_1$; when $\kappa=\aleph_\epsilon$ and  the
transcendence degree of $A$ is infinite (thus the case of finite
transcendence degree of $A$ is left open).\\[0.05in]
These results are not surprising, given that any complete theory of
pseudo-finite fields has the independence property. However, the proofs
do use some properties which are specific to pseudo-finite fields, so it
is not clear that the results would hold in the general case of theories
with IP. The question arose
during the study of the existence (and uniqueness) of certain
strengthenings of the notion of difference closure of difference fields
of characteristic $0$. In \cite{Ch}, we show that if $K$ is an
algebraically closed difference field of characteristic $0$, and
$\kappa$ an uncountable cardinal or $\aleph_\varepsilon$, and if
$\fix(\si)(K)$ is a  $\kappa$-saturated pseudo-finite field,
then $\kappa$-prime models of ACFA (the theory of existentially closed
difference fields) over $K$ exist and are unique up to
$K$-isomorphism. The question then arises of whether the hypothesis on
the fixed field of $K$ is necessary. This note shows that it
is, under natural assumptions. \\[0.05in]
The paper is organised as follows. Section 1 recalls well-known facts
about fields, Section 2 gives the results on the non-existence of prime
models, and Section 3 those on the non-existence of $\kappa$-prime models.

\section{Preliminaries}

\para {\bf Convention and notation}. Unless otherwise mentioned,
all fields will be subfields of a large algebraically closed field. If $K$ is a field, then $K^s$ denotes the
separable closure of $K$, $K^{alg}$ its algebraic closure, and $G(K)$
its absolute Galois group $\gal(K^s/K)$. If $L$ is an extension of
the field $K$, and $\si\in \aut(L/K)$, then $\fix(\si)$ will denote the
subfield of $L$ consisting of elements fixed by $\si$. If $\si\in G(K)$, then  $\langle \si\rangle$ denotes the
topological closure inside $G(K)$ of the group generated by $\sigma$.

\para {\bf Classical algebraic results on fields}. \label{reg} (See chapter 3 of Lang's book \cite{La}) Let $K\subset L$ be fields. Recall that
$L$ is {\em regular} over $K$ if it is linearly disjoint from $K^{alg}$
over $K$. If $K$ is {\em perfect} (i.e., of characteristic $0$, or if of
characteristic $p>0$, closed under $p$-th roots), then this is
equivalent to $L\cap K^s=K$. The {\em perfect hull} of $K$ is $K$ if
char$(K)=0$, and the closure of $K$ under $p$-th roots if
char$(K)=p>0$.
The field $L$ is {\em separable} over $K$
if it is linearly  disjoint from the perfect hull of $K$ over
$K$. Finally, if $L$ is separable over $K$, then $L\cap K^s=K$ implies
that $L$ is regular over $K$. \\
Recall also that a polynomial $f\in K[\bar X]$ is called {\em absolutely
  irreducible} if it is irreducible in $K^{alg}[\bar X]$. This
corresponds to the field ${\rm Frac}(K[\bar X]/(f))$ being a regular
extension of $K$.

\para{\bf The Haar measure}. \label{haar} Recall that if $K$ is a field, then $G(K)$
can be endowed uniquely with a measure $\mu$ on the $\si$-algebra
generated by open subsets of $G(K)$, which satisfies $\mu(G)=1$, and is
stable under translation. This measure is called the {\em Haar
  measure}. If $L$ is a finite separable extension of $K$, then
$\mu(G(L))=[L:K]\inv$. Furthermore, assume that $L_i$, $i<\omega$, is a
family of linearly disjoint algebraic extensions of $K$ and $A_i$ a non-empty set of left-cosets of
$G(L_i)$ in $G(K)$. If
$\sum_i[L_i:K]\inv =\infty$, then $\mu(\bigcup_i A_i)=1$ (Lemma~18.5.2 in
\cite{FJ}).

\para {\bf Review on Hilbertian fields and their
  properties}. \label{hilb1} All
references are to the book of Fried and Jarden,  \cite{FJ}.
\begin{enumerate}
\item Recall that a field $K$ is {\em Hilbertian} if whenever $f\in
  K[T,X]$ ($(T,X)$ a tuple of indeterminates, $|X|=1$) is separable in $X$ and irreducible over $K(T)[X]$, then there are infinitely
  many tuples $a$ in $K$ such that $f(a,X)$ is  irreducible over $K$.
There are many equivalent statements of this property, and in particular
if it is satisfied for $|T|=1$, then it is satisfied for tuples
$T$ of arbitrary length (Proposition 13.2.2).
\item Examples of Hilbertian fields include $\rat$ and any finitely generated
 infinite field. Function fields are Hilbertian, and if $K$ is
 Hilbertian, then so is any finite algebraic extension of $K$. An
 infinite separably algebraic extension $L$ of a Hilbertian field $K$ is  not
 necessarily Hilbertian, but any finite proper separable extension of
 $L$ which is
 not contained in the Galois hull of $L$ (over $K$) is Hilbertian (Theorem 13.9.4).
\item Let $M_1,M_2$ be Galois extensions of the Hilbertian field $K$, and
$M$ a subfield of $M_1M_2$ containing $K$ and such that $M\not\subset M_i$ for $i=1,2$. Then $M$
is Hilbertian (Theorem 13.8.3).
\item In order to state some properties of Hilbertian fields, it is convenient to define, for $K$ a field and irreducible
polynomials $f_1,\ldots, f_m\in
K[T][X]$ which are separable in $X$, and non-zero $g\in K[T]$, the {\em separable Hilbert set} $H_K(f_1,\ldots,f_m; g)$ as the set of
$a\in K$ such that $g(a)\neq 0$ and $f_1(a,X),\ldots,f_m(a,X)$ are
irreducible over $K$.
\item Every separable Hilbert subset of $K^r$ contains one of the form $H_K(f)$, with
$f$ monic irreducible and separable (Lemma 12.1.6). Hence if $K$ is Hilbertian then
  every separable Hilbert set is infinite.
\item Let  $L$ be a finite separable extension
of $K$. Then every separable Hilbert subset of $L$ contains a separable
Hilbert
subset of $K$ of the forme $H_K(f)$ (Lemma 12.2.2).
\item Let $K$ be a Hilbertian field, $f(T,X)\in K[T,X]$ irreducible and
separable in $X$, and $G$ the Galois group of the Galois extension of
$K(T)$ generated by the roots of $f(T,X)=0$. Then there is a separable
Hilbert set $H\subseteq K^r$ such that if $a\in H$, then the Galois group
of the extension generated by the roots of $f(a,X)=0$ is isomorphic to
$G$. In particular, $f(a,X)$ is irreducible (Proposition 16.1.5).

\end{enumerate}

\begin{fact} \label{rem1} Some easy observations and reminders about fields.
\begin{enumerate}
\item Let $B$ be a {\em primary}\footnote{i.e.,  $B\cap A^{s}=A$.} extension of the field $A$, and $\si\in
  G(A)$. Then $\si$ lifts to some $\si'\in G(B)$. Indeed, $\si$ has an
  obvious extension to $A^s\otimes_AB$ given by $\si'=\si\otimes id$;  by
  primarity of $B/A$, $A^s\otimes_AB $ is a domain, and
  is isomorphic to $A^sB$. This automorphism $\si'$ of $A^sB$ extends to an
  automorphism of $B^s$ which is the identity on $B$. Recall that if $B$
  is a regular extension of $A$ then it is primary.
\item Let $K$ be a field, $\si\in G(K)$. Then $\langle
  \si\rangle\simeq\hat\zee$ if and only if $\langle
  \si\rangle$ has a quotient isomorphic to $\zee/4\zee$, and quotients
  isomorphic to $\zee/p\zee$ for every odd prime $p$. The necessity is
  clear, the sufficiency follows from the fact that the only possible
  order of a torsion
  element of the absolute Galois group of a field is $2$ (and then the
  field is of characteristic $0$ and does not contain $\sqrt{-1}$) and that
  $\langle\si\rangle$ is the direct product of its Sylow subgroups. When
  char$(K)$ is positive, it suffices that $G(K)$ has a quotient
  isomorphic to $\zee/p\zee$ for every prime $p$.

\item Recall that by Theorem 11.2.3 of \cite{FJ}, if  $L$
  an algebraic extension of a field $K$, and every absolutely irreducible
  $f(X,Y)\in K[X,Y]$ has a zero in $L$, then $L$ is PAC.

\item Let $K$ be a field, $f(X,Y)\in K[X,Y]$ an absolutely irreducible
  polynomial. If $f(X,Y)$ is not separable as a polynomial in $X$, then
  it is separable as a polynomial in $Y$. Indeed, otherwise it would not
  stay irreducible over the perfect hull of $K$.
\item (Kummer theory). Let $K$ be a field of characteristic not $2$, let $t$ be
  transcendental over $K$, and $a_1,\ldots,a_n$ distinct elements of
  $K$. Then the fields $K(t)(\sqrt{t+a_i})$ are linearly
  disjoint over $K(t)$, and they are proper Galois extensions of
  $K(t)$. Moreover, the field $K(t)(\sqrt{t+a_i}\mid 1\leq i\leq n)$ is
  a regular extension of $K$. The general phenomenon is as follows: let
  $L$ be an extension of
  $K(t)$ generated by square roots of polynomials $f_i(t)$,
  $i=1,\ldots,n$, and assume that the elements $f_i(t)$ are
  multiplicatively independent modulo the multiplicative subgroup
  $K^\times {K(t)^\times}^2$ of $K(t)^\times$; then $L$ is a regular extension of
  $K$, and $\gal(L/K(t))\simeq (\zee/2\zee)^n$.
\item  (Artin-Schreier theory) Let $K$ be a field of characteristic $2$, and $a_1,\ldots,a_n\in
  K$ be $\ffi_2$-linearly independent. Let $\alpha_i$ be a root of
  $X^2+X+a_it=0$ for $i=1,\ldots,n$. Then the fields $K(t)(\alpha_i)$
  are linearly disjoint over $K(t)$, and are proper Galois extensions of
  $K(t)$. Moreover, the field $K(t,\alpha_1,\ldots,\alpha_n)$ is a
  regular extension of $K$. The general phenomenon is as follows: let
  $L$ be an extension of
  $K(t)$ generated by solutions of $X^2+X+f_i(t)=0$, $i=1,\ldots,n$, where the $f_i(t)$
  are elements of $K(t)$, which are $\ffi_2$-linearly independent modulo
  the additive subgroup $K+\{f(t)^2-f(t)\mid f(t)\in K(t)\}$ of $K(t)$;
  then $L$ is a regular extension of
  $K(t)$, and $\gal(L/K(t))\simeq (\zee/2\zee)^n$.

\item (Linear disjointness). Recall that  if $M\subset N$ are
fields, and $L$ is a Galois extension of $M$, then $N$ and $L$ are linearly
disjoint over $N\cap L$. The same holds if $N$ is perfect and $L$ is
the perfect hull of a Galois extension of $M$, because $L$ will then be
a Galois
extension of the perfect field $L\cap N$.
This remark will be constantly used.
\item (Theorem III.3 in \cite{La}) Let $K/k$ be a regular extension of
  fields, let the field $L$ contain
  $k$, and assume that $K$ and $L$ are free over $k$. Then $K$ and $L$
  are linearly disjoint over $k$.

\end{enumerate}
\end{fact}

\begin{lem}\label{lem0} Let $K\subset L_1,L_2$ be three algebraically closed fields, with
  $L_1$ and $L_2$ linearly disjoint over $K$, and consider the field
  composite $L_1L_2$. Let $u\in L_1\setminus K$,
  $v\neq w\in L_2\setminus K$.
\begin{enumerate}
\item If char$(K)\neq 2$, then
  $[L_1L_2(\sqrt{u+v},\sqrt{u+w}):L_1L_2]=4$. %$\sqrt{u+v}\notin L_1L_2$.
\item Assume char$(K)=2$, that $v+w\in L_2\setminus K$, and let $c$
  be a root of $X^2+X+uv$, $d$ a root of $X^2+X+uw$. Then
  $[L_1L_2(c,d):L_1L_2]=4$.
\end{enumerate}
\end{lem}

\prf (I thank  Olivier Benoist for this elegant proof.) \\
(1) In characteristic $\neq 2$, it suffices to prove that neither
$c=\sqrt{u+v}$, nor $\sqrt{u+v}\sqrt{u+w}$ is in $L_1L_2$; and in
characteristic $2$ that neither $c$ nor $c+d$ is in $L_1L_2$.

Let us first do the case of odd characteristic. Assume by way of
contradition that both $\sqrt{u+v}$ and $\sqrt{u+v}\sqrt{u+w}$ belong to
$L_1L_2$. Then there are finite tuples $u_1\in L_1$ and $v_1\in L_2$
such that $u+v$ and $(u+v)(u+w)$ have square roots in
$K(u,v,w,u_1,v_1)$. As $L_1$ and $L_2$ are free over $K$, $\tp_{\rm
  ACF}(v,w, v_1/L_1)$ does not fork over
$K$ (in the sense of the theory ACF of algebraically closed fields), and
therefore is finitely satisfiable in $K$. In particular, there are
infinitely many pairs $b_1\neq c_1$ in $K$ such that both $u+b_1$ and
$(u+b_1)(u+c_1)$ have a square root in the field $K(u,u_1)$. But this is
impossible: as we saw above in \ref{rem1}(5), the extensions $K(\sqrt{u+b})$, $b\in K$,
are linearly disjoint over $K(u)$, and therefore $K(u,u_1)$ contains at
most finitely many of them, since it is finitely generated over
$K(u)$.\\
(2) Same proof: assume that both $c$ and $c+d$ are in $L_1L_2$, and let
$u_1\in L_1$ and $v_1\in L_2$ be finite tuples such that $c,c+d\in
L_1L_2$. Then for infinitely many pairs $(e,f)$ in $K$ which are
$\ffi_2$-independent, we would have that both $X^2+X+ue$ and
$X^2+X+u(e+f)$ have a solution in $K(u,u_1)$ which is impossible by \ref{rem1}(6).

%% If $\alpha=\sqrt{u+v}\in L_1L_2$, then there are finite tuples
%% $u_1\in L_1$ and $v_1\in L_2$ such that $\alpha\in K(u,u_1,v,v_1)$. As
%% $L_1$ and $L_2$ are free over $K$, $\tp_{\rm ACF}(v,v_1/L_1)$ does not fork over
%% $K$ (in the sense of the theory ACF of algebraically closed fields), and
%% therefore is finitely satisfiable in $K$. In particular there are
%% infinitely many elements
%% $b\in K$ such that $u+b$ has a square root in $K(u,u_1)$.
%% But this is
%% impossible: as we saw above in \ref{rem1}(5), the extensions $K(\sqrt{u+b})$, $b\in K$,
%% are linearly disjoint over $K(u)$, and therefore $K(u,u_1)$ contains at
%% most finitely many of them, since it is finitely generated over $K(u)$. \\
%% (2) Same proof: If $\alpha\in L_1L_2$, then $\alpha\in K(u,u_1,v,v_1)$
%% for some $u_1\in L_1$ and $v_1\in L_2$. For all but finitely many $b\in
%% K$, the equation $X^2+X+bu$ would have a solution in $K(u,u_1)$, but
%% this is impossible by \ref{rem1}(6).

\begin{cor}\label{cor0} Let $L_1$ and $L_2$ be regular extensions of the field $K$,
  which are linearly disjoint over $K$. Let $u\in L_1\setminus K$,
  $v\neq w\in
  L_2\setminus K$.
\begin{enumerate}
\item If char$(K)\neq 2$, then $[L_1L_2(\sqrt{u+v},\sqrt{u+w}):L_1L_2]=4$
  and $L_1L_2(\sqrt{u+v},\sqrt{u+w})\cap
  L_1^{alg}L_2^{alg}=L_1L_2$. Hence $L_1L_2(\sqrt{u+v},\sqrt{u+w}))$ is a regular
  extension of both $L_1$ and $L_2$.
\item Assume char$(K)=2$, that $v$, $w$, $v+w$ is in $K$, and let $c$ be
  a root of $X^2+X+uv$, $d$ be a root of $X^2+X+uw$. Then
  $[L_1L_2(c,d):L_1L_2]=4$, $L_1L_2(c,d)\cap L_1^{alg}L_2^{alg}=L_1L_2$, and $L_1L_2(c,d)$ is a regular extension of both $L_1$ and  $L_2$.
\end{enumerate}
\end{cor}

\prf (1) Our assumption implies that $L_1$
and $L_2$ are free over $K$, and therefore that their algebraic closures
$L_1^{alg}$ and $L_2^{alg}$
are linearly disjoint over $K^{alg}$. By Lemma \ref{lem0},
$L_1^{alg}L_2^{alg}(\sqrt{u+v},\sqrt{u+w})$ has maximal degree $4$ over
$L_1^{alg}L_2^{alg}$, so $L_1L_2(\sqrt{u+v},\sqrt{u+w})$ has also degree
$4$ over $L_1$ and $L_2$, and is therefore a regular extension of bothe
$L_1$ and $L_2$.
%
%% \sqrt{u+v}\notin L_1^{alg}L_2^{alg}$. I.e.,
%% $L_1L_2(\sqrt{u+v})\cap L_1^{alg}L_2^{alg}=L_1L_2$, so that
%% $L_1L_2(\sqrt{u+v})$ is a regular extension of both $L_1$ and $L_2$.
Same proof for (2).

\begin{lem} \label{lem1} (Folklore) Let $G$ be a finite abelian group, $F$ a
  field, $t$ an
indeterminate, and assume that $F$ has only finitely many  Galois
extensions with Galois group isomorphic to a quotient of $G$. Then there is a sequence $L_i$, $i\in\omega$, of linearly
disjoint Galois extensions of $F(t)$ with Galois group isomorphic to
$G$, and the field composite of which is a regular extension of $F$.

\end{lem}

\prf Let $M$ be the composite of the finitely many abelian Galois
extensions of $F$ with Galois group isomorphic to a quotient of $G$. Let $u$ be a new indeterminate. By
Proposition 16.3.5 of \cite{FJ}, letting $K=F(t)$, the field $K(u)$ has a Galois
extension $L$ which is regular over $K$, and with Galois group $G$. Let
$\alpha$ be a generator of $L$ over $K(u)$, and $f(u,X)\in K(u)[X]$ its
minimal polynomial over $K(u)$. As $L$ is regular over $K$, $f(u,X)$ is
irreducible over $M(t,u)$. \\
Observe that if $L'$ is a Galois extension of $K$ with Galois group
$G$, and if $L'\cap M=F$, then $L'$ is regular over $F$. Indeed,
$L'\cap F^s$ is a Galois extension of $F$, with Galois group isomorphic
to a
quotient of $G$, and therefore is contained in $M$. Our assumption
therefore
implies that $L'\cap F^s=F$. Furthermore, $L'$ is separable over $F$, hence
regular over $F$. \\
As $K$ is Hilbertian, by Property \ref{hilb1}(6) and (7) there is $a\in K=F(t)$ such that $f(a,X)$ is
irreducible over $M(t)$, and such that the field $L_0$ generated over $F(t)$
by a root of $f(a,X)$ is Galois with Galois group isomorphic to
$G$. Then $L_0\cap M=F$, and by
the discussion in the previous paragraph, $L_0$ is a  Galois
extension of $F(t)$ which is regular over $F$. \\
Replacing $M(t)$ by $ML_0$, we construct in the same fashion a Galois
extension $L_1$ of $K$, with Galois group isomorphic to $G$, and which
is linearly disjoint from $ML_0$ over $K$. We iterate the construction
and build by induction a sequence $L_i$, $i\in\nat$, of Galois
extensions $L_i$ of $K$ with Galois group isomorphic to $G$, and such
that for every $i$, $L_i$ is linearly disjoint from $ML_0\cdots L_{i-1}$
over $K$. In particular, the field composite of all $L_i$'s is a regular
extension of $F$.

%% \begin{lem}\label{lem2} Let $M$ and $C$ be perfect fields, regular and linearly
%%   disjoint over their intersection $M\cap C$. Assume that $G(M)$ and
%%   $G(C)$ are procyclic, and that $a\in (MC)^{alg}$ is such that $MC(a)$
%%   is regular over $M$ and over $C$.
%% \begin{enumerate}
%%  \item Then there is $\si\in G(MC)$ such that $\fix(\si)$ is regular
%%    over $M$ and over $C$ and contains $a$.
%% \item
%%   Furthermore, if $MC(a)\cap M^{alg}C^{alg}=MC$, and $a'$ is an
%%   $MC$-conjugate of $a$, then there is $\si\in G(MC)$, such that
%%   $\fix(\si)$ is a regular extension of $M$ and of $C$, and such that
%%   $\si(a)=a'$.

%%   \end{enumerate}
%% \end{lem}

%% \prf

\para {\bf Review on pseudo-finite fields and their
  properties}. \label{psf1} Recall
that the theory of pseudo-finite fields is axiomatised by the following
properties: the field is {\em PAC} (every absolutely irreducible variety
defined over the field has a rational point); the absolute Galois group
is isomorphic to $\hat\zee$ ($=\lim_{\leftarrow}\zee/n\zee$); if the
characteristic is $p>0$, then the field is {\em perfect} (closed under
$p$-th roots). We will mainly use the following five results:
\begin{enumerate}\item  Let $F_1$ and $F_2$ be two pseudo-finite fields containing a common
subfield $E$. Then
$$F_1\equiv_E F_2 \hbox{ if and only if there is an }E\hbox{-isomorphism
} F_1\cap
E^{alg}\to F_2\cap E^{alg}.$$

  \item Let $L$ be a relatively algebraically closed subfield of the perfect field
    $E$ and of the $|E|^+$-saturated
    pseudo-finite field $F$. Assume that $G(E)$ is procyclic. Then there
    is an $L$-embedding $\Phi$ of $E$ into $F$ such that
    $F/\Phi(E)$ is regular.
    \item If $E$ is a perfect field with procyclic absolute Galois group,
  then it has a { regular} extension  $F$ which is pseudo-finite.
\end{enumerate}
(1) is a special case of 20.4.2 in \cite{FJ}. \\
(2) follows from the
Embedding Lemma (20.2.2 and 20.2.4 in \cite{FJ}) with $\Phi_0=id$: the
restrictions maps ${\rm res}_{F/L}:
G(F)\to
G(L)$ and ${\rm res}_{E/L}: G(E)\to G(L)$ are onto, and because $G(E)$ is procyclic and
$G(F)$ is free, there is an onto map $\varphi: G(F)\to G(E)$ such that
${\rm res}_{E/L}\varphi={\rm res}_{F/L}$. The lemma then gives the map $\Phi$, and
because $\varphi$ is onto and $E$ is perfect, the extension $F/\Phi(E)$ is regular.\\
(3) is folklore, but I was
not able to find an explicit statement of it: when $E$ is a subfield of
the algebraic closure of the prime field $k$, this is given by
Propositions 7 and 7' of \cite{Ax}. In the general case, $E$ is a
regular extension of $L:=k^{alg}\cap E$, and $L$ has procyclic
Galois group and is perfect. By the above, there is some pseudo-finite
field $F$ containing
$L$, which is regular over $L$, and we may assume it is sufficiently
saturated. Because $G(E)$ is procyclic, there is an onto map $\varphi: G(F)\to G(E)$ such that
${\rm res}_{E/L}\varphi={\rm res}_{F/L}$, and we conclude as in (2). \\[0.1in]
These three results have several consequences. For instance, if $E\subset F_1$
is relatively algebraically closed in the pseudo-finite field $F_1$,
then the theory Psf
together with the quantifier-free diagramme of $E$ is complete
 (in the language $\call(E)$ of rings augmented by
constant symbols for the elements of $E$). \\
In particular, if $a\in F_1$ is transcendental over $E$, then $\tp(a/E)$
is entirely axiomatised by the collection of $\call(E)$-formulas
expressing that it is transcendental over $E$, as well as, for each
finite Galois extension $L$ of $E(a)$, a formula which describes the
isomorphism type over $E(a)$ of $L\cap F_1$. So this formula will say
which polynomials $f(a,X)\in E[a,X]$ have a solution in $F_1$ and which
do not. By (3) above, note that any  subfield $K$ of $L$ which is a
regular extension of $E$ and with $\gal(L/K)$ cyclic can appear as
$L\cap F$ for some model $F$ of  $T(E)$ which
contains $a$.
\begin{enumerate}
\item[(4)] Hence, if $F_1$ is a pseudo-finite field containing $E$ and
  regular over $E$, and $a\in F_1$ is transcendental over $E$, then
  $\tp(a/E)$ is not isolated. This follows easily from the description of
  types, and  because $E(a)$ has infinitely many linearly disjoint
  extensions $L_i$ ($i\in\nat$), the composite $L$ of which  is regular
  over $E$ (see Lemma \ref{lem1}). Indeed the
  type of $a$ is axiomatized by saying that $a$ is trancendental over
  $E$, and by saying which polynomials $f(a,X)\in E[a,X]$ have a root
  in $F_1$ and which have not. In particular, any $\call(E)$-formula
  $\varphi(x)$ will only give
  information about $F_1\cap L_0\cdots L_n$ for some  $n$, and say
  nothing about $F_1\cap L_{n+1}$, and whether it equals $L_{n+1}$ or
  not.

\item[(5)] If $F_1$ is pseudo-finite and $E\subset F_1$, then $\acl(E)=E^{alg}\cap F_1$, see Proposition~4.5 in
  \cite{CP}.
\end{enumerate}

\section{Non-existence of prime models}

\para {\bf Setting}. \label{setting} Let $T$ be a complete theory of pseudo-finite
fields, $\ffi$ a model of
$T$, and $A\subset \ffi$, $T(A)$ the $\call(A)$-theory of $\ffi$ $(\call$ the
language of rings $\{+,-,\cdot,0,1\}$). We want to
show that unless $\acl(A)$ is a pseudo-finite field, then $T(A)$ has no
prime model. As $T(A)$
describes the $A$-isomorphism type of $\acl(A)=A^{alg}\cap \ffi$ over $A$,
without loss of generality, we will assume that $A^{alg}\cap
\ffi=A$. Note that $A$ is perfect, $G(A)$ is procyclic, and we
will fix a topological generator $\si$ of $G(A)=\gal(A^s/A)$.

\begin{notation} Let $A$ be a field, $F$ a regular field extension of
  $A$, and $t\in A$. We
  denote by $\cals(t,F)$ the set
$$\cals(t,F)=\begin{cases} \{a\in A\mid \sqrt{t+a}\in F\} \hbox{ if }
  {\rm char}(A)\neq 2,\\
\{a\in A\setminus\{0\} \mid F\models \exists y\, y^2+y=at\} \hbox{ if }{\rm char}(A)=
2.
\end{cases}$$
\end{notation}

\begin{rem} \label{rem2} Observe that if $F\subseteq F'$,
  $F^{alg}\cap F'=F$ and $t\in F$, then
  $\cals(t,F)=\cals(t,F')$.
\end{rem}

\begin{prop} \label{prop1} Let $T$ and $A$ be as above, with $A$ not
  pseudo-finite. Then $T(A)$ has a model $F_0$ of transcendence degree
  $1$ over $A$. Furthermore:
\begin{itemize}
\item[(1)] Assume that $A$ is countable, let $t$ be transcendental over
  $A$, and let $\tilde\si$ be a
  lifting of $\si$ to $G(A(t))$. Then for almost all $\tau\in
  G(A^s(t))$, the perfect closure of the subfield of $A(t)^s$ fixed by $\tilde\si\tau$ is a
  model of $T(A)$.
\item[(2)] Assume that $|A|=\kappa\geq \aleph_0$. When char$(A)\neq 2$,
  we choose some $X\subset A\setminus \{0\}$. If char$(A)=2$, we
  fix a basis $Z$ of the $\ffi_2$-vector space $A$ with $1\in Z$,
  and take $X\subset Z$.
  Then there is a model $F_X$ of $T(A)$ which has
  transcendence degree $1$ over $A$, and is such that for some $t\in
  F_X\setminus A$,
$$\begin{cases}
\cals(t,F_X)=X &\hbox{ when char}(A)\neq 2,\\
\cals(t,F_X)\cap Z=X & \hbox{ when char}(A)= 2.
\end{cases}
$$

\end{itemize}
\end{prop}

\prf A model  of $T(A)$ is a regular extension of $A$, with
absolute Galois group isomorphic to $\hat\zee$, and which is PAC and
perfect. For
both items we will construct the model as an algebraic extension of $A(t)$: we
will first work inside $A(t)^s$, then take the perfect closure.
Recall that by \ref{rem1}(3), for the PAC condition, it suffices to build a regular
extension of $A$ contained in $A(t)^s$, and in which  every
absolutely irreducible plane curve defined over $A(t)$ has a point. Then
its perfect closure will be pseudo-finite. We first show (1). We
will show the following:
\begin{itemize}[noitemsep]
\item[(i)] if $f(X,Y)\in A(t)[X,Y]$ is absolutely
irreducible, then for almost all $\tau\in G(A^s(t))$ (in the sense of
the Haar measure $\mu$ on $G(A^s(t))$),
$\fix(\tilde\si\tau)$ contains a solution of $f(X,Y)=0$.
\item[(ii)] for almost all $\tau\in G(A^s(t))$, for every $n\geq 2$,
$\langle \tilde
\si\tau\rangle$ has a quotient isomorphic to $\zee/n\zee$.
\end{itemize}
Towards (i), let $f(X,Y)\in A(t)[X,Y]$ be absolutely irreducible; by
Fact \ref{rem1}(4) we may assume that $f$ is separable in $Y$, and we
let $m$ be the degree of $f$ in $Y$. Let $B$ be the subfield of $A^s$
fixed by ${\si}^{m!}$. As $A(t)$ is Hilbertian, as in the proof
of Lemma \ref{lem1} (using Property \ref{hilb1}(6)), we  build
inductively a sequence $L_i$, $i<\omega$, of finite separable extensions of
$A(t)$, and of elements $a_i\in A(t)$, such that: \\
-- the polynomial $f(a_i,Y)$ is irreducible over   $BL_0\cdots L_{i-1}$
for all $i$ (over $B(t)$ if $i=0$);\\
-- $L_i=A(t,b_i)$ where $f(a_i,b_i)=0$.\\[0.05in]
(For more details  one may look at Theorem 18.6.1 in \cite{FJ}.) Note that because
$[L_i:A(t)]\leq m$, it follows that $L_i$ is linearly disjoint from
$A^sL_0\cdots L_{i-1}$ over $A(t)$ for every $i$, and therefore that the field
composite $L$ of all $L_i$'s is a regular extension of $A$. By Fact
\ref{rem1}(1),
$\si$ extends to some $\si'\in G(L)$. Then, for every $\tau\in \bigcup
_i G(A^sL_i)$, $\fix(\si'\tau)$ contains a solution of
$f(X,Y)=0$. Hence, for every $\tau\in ({\tilde \si}\inv \si')(\bigcup_i
G(A^sL_i))$, $\fix(\tilde \si \tau)$ contains a solution of
$f(X,Y)=0$. By \ref{haar}, $\mu(\bigcup_i G(A^sL_i))=1$,
and so does its translate by ${\tilde \si}\inv \si'$. This shows (i).\\[0.05in]
(ii) is proved in the same fashion, using \ref{rem1}(2). Let $n$ be a
prime or $4$, and use Lemma \ref{lem1} to find a sequence $(L_i)_{i<\omega}$ of linearly
disjoint Galois extensions of $A(t)$, with $\gal(L_i/A(t))\simeq
\zee/n\zee$, and such that the field composite $L$ of all $L_i$'s is a
regular extension of $A$.
As in (i), the
set of $\tau\in G(A^s(t))$ such that for some $i$, $\tau\rest_{L_i}$
generates $\gal(L_i/A(t))$, has measure $1$, and therefore so does its
translate (on the left) by ${\tilde\si}\inv\si'$. This proves (ii). \\

\noindent
A countable intersection of sets of Haar measure $1$ has measure $1$,
and therefore the set of $\tau\in G(A^s(t))$ such that
\begin{center} every absolutely irreducible $f(X,Y)$ has a solution in
$\fix(\tilde\si\tau)$, and $\langle\tilde\si\tau\rangle \simeq
\hat\zee$ \end{center}
has measure $1$. For any such $\tau$, the field $\fix(\tilde\si\tau)$ is
therefore PAC, with absolute Galois group isomorphic to $\hat \zee$, and
 its perfect closure is our desired pseudo-finite field. \\[0.05in]
(2) There are four cases to consider, depending on the characteristic,
and whether $A$ has an algebraic extension of degree $2$ or
not. Let $t$ be an indeterminate over $A$.\\[0.05in]
Case 1: char$(A)\neq 2$ and  $A^2\neq A$: \\
Let $c\in
A\setminus A^2$, and consider the Galois extension $L_0$ of $A(t)$ defined
as the field composite of all $A(t,\sqrt{t+a})$ for $a\in X$, and
$A(t,\sqrt{ct+ca})$ for $a\in A\setminus X$. If $B\supseteq L_0$ is
regular over $A$, then $\cals(t,B)=X$: indeed, $c(t+a)\in B^2$, $c\notin
B^2$ imply $(t+a)\notin B^2$.  \\[0.05in]
Case 2: char$(A)\neq 2$ and $A^2=A$:\\
We let $L_0$ be the field composite
of all $A(t,\sqrt{t+a})$ for $a\in X$, and all $A(t,\sqrt{t^2+at})$ for
$0\neq a\in A\setminus X$. If $B\supseteq L_0$ is such that
$t\notin B^2$, then $\cals(t,B)=X$. \\[0.05in]
Case 3: char$(A)=2$, and $A$ has an extension of degree
$2$, say $Y^2+Y+c=0$ has no solution in $A$: \\
Let $L_0$ be the field obtained by
adjoining to $A(t)$ a solution of $Y^2+Y+at=0$ if $a\in X$, and a solution of
$Y^2+Y+at+c=0$ if $a\in Z\setminus X$. Then if $B\supseteq L_0$ is a regular
extension of $A$, we have $\cals(t,B)\cap Z=X$. \\[0.05in]
Case 4: char$(A)=2$, and  $A$ is closed under
Artin-Schreier extensions:\\
Let $\alpha$ satisfy $Y^2+Y+t^3=0$, and let $L_0$
be the Galois extension of $A(t)$ obtained by adjoining a solution of
$Y^2+Y+at=0$ if $a\in X$, and $Y^2+Y+at+t^3=0$ if $a\in Z\setminus X$. Again, if
$B\supseteq L_0$ does not contain $\alpha$, then $\cals(t,B)\cap Z=X$. \\[0.05in]
Note that in all four cases, $L_0$ is regular over $A$ (by Facts \ref{rem1}(5)
and (6)), and is Hilbertian (by Property \ref{hilb1}(3)). It therefore
suffices to construct an algebraic
extension of $L_0$ which is regular over $A$, does not contain the forbidden elements $t^{1/2}$ or $\alpha$ when $A$ has no
proper algebraic extension of degree $2$, and is pseudo-finite. To
do the latter, we will  construct inside $L_0^s$ a PAC field which
contains $L_0$, and with
Galois group isomorphic to $\hat \zee$. We first take care of the
Galois group. To do that, we will find some Galois extension $L$
of $A(t)$, which is linearly disjoint from $A^sL_0$ over
$A(t)$, and  such
that $\gal(L/A(t))\simeq \hat\zee$. Let
$Q$ be the set of $n$ which are prime numbers or $4$  and such that
$G(A)$ does not have a quotient isomorphic to $\zee/n\zee$. Note that
$2\notin Q$ and $4\in Q$ implies that the characteristic is $0$,
$i\notin A$, and $A(i)$ contains all $2^n$-th roots of unity. \\[0.05in]
For each
odd $n\in Q$, using Lemma \ref{lem1} we find a Galois extension $L_n$ of $A(t)$ with Galois group
isomorphic to $\zee/n\zee$ and which is regular over $A$. Note
  that automatically, the field composite of all $L_n$ (with $n$ odd in $Q$) will be linearly
disjoint from $L_0$ over $A(t)$.  When $n$ is
$2$ or $4$ we will need to be a little more careful. \\[0.05in]
Case 3 is vacuous, as is
Case 1 when char$(A)\neq 0$.
In case 2, $A$ contains $\sqrt{-1}$, and we let $L_2=L_4=A(t^{1/4})$. Then $L_2$ is linearly
disjoint from $A^sL_0$ over $A(t)$, with Galois group $\zee/4\zee$.  In
case 4, we let $L_2=L_4=A(t)(\alpha)$; it is linearly disjoint from
$A^sL_0$ over $A(t)$.\\
We are left with Case 1, char$(A)=0$, $4\in Q$, $2\notin Q$, and
  therefore $\sqrt{-1}\notin A$. This case %% of $2\notin Q$ and $4\in Q$, so that
%% char$(A)=0$,
is more delicate, and we proceed as follows
(it is a particular case of the construction
given in Lemma~16.3.1 of \cite{FJ}). We fix a square root $i$ of $-1$; then
$\si(i)=-i$ ($\si$ a generator of $G(A)$). Consider the element $1+it$,
and let $a \in A(t)^s$ satisfy $a^4=1+it$. Such an element $a$ can be
found in $A(i)[[t]]$ (by Hensel's lemma), and we may therefore lift
$\si\rest_{A(t,i)}$ to an element $\si_1\in \aut(A(t,i,a)/A(t))$ with $\si_1^2=id$. Let
$b=a\si_1(a)^3$, and
note
that $$b^4=(1+it)(1-it)^3=(1+t^2)(1-it)^2,$$ and that $1+t^2$ has no
square root in $L_0(i)$ (as $1+t^2=(1+it)(1-it)$ is relatively prime to
all $(1+ta)$ with $a\in A$ and by Kummer theory -- Fact \ref{rem1}(5)).
By definition of $L_0$, we have
$$L_0(i,\sqrt{1+it},\sqrt{1-it})= A(t,i)(\sqrt{t+a}\mid a\in X)(\sqrt{i(t+a)}\mid a\in
A\setminus X)(\sqrt{1+it},\sqrt{1-it}).$$
Note that we are taking square roots of polynomials of degree $1$ over
$A(i)$, and that they are all relatively prime, so that by Fact
\ref{rem1}(5), this field is a regular extension of $A(i)$, and moreover
$\sqrt{1+t^2}\notin L_0(i)$.

Hence $[A(t,i,b):A(t,i)]=4$, and $[A(t,i,b):A(t)]=8$. Define $\omega\in \gal(A(t,i,b)/A(t,i))$ by $\omega(b)=ib$.  We now compute $\omega\si_1$ and $\si_1\omega$ on $i$ and
on $b$. We have:
$$\omega\si_1(i)=\omega(-i)=-i, \qquad \si_1\omega(i)=\si_1(i)=-i, \qquad \si_1\omega(b)=\si_1(ib)=-i\si_1(b)$$
and one computes
$$\omega\si_1(b)=\omega(\si_1(a)a^3)=\omega(b^3 \si_1(a)^{-8})=-ib^3
\si_1(a)^{-8} =-i \si_1(b).$$
(Here we use that $\si_1$ is an involution, that $\si_1(a)^8\in A(t,i)$ is
fixed by $\omega$). So $\si_1$ and $\omega$ commute, and  $\gal(A(t,i,b)/A(t))$ is the direct product of
the subgroups generated by $\si_1$ and by $\omega$. We take $L_4$ to be
the subfield of $A(t,i,b)$ fixed by $\si_1$. It is regular over
$A$, with Galois group over $A(t)$ isomorphic to $\zee/4\zee$. One
computes that $\sqrt{1+t^2}\in A(t,i,b)$ is fixed by $\si_1$ (since
$1+t^2=(a^2\si_1(a)^2)^2$), and therefore $L_4$ is linearly
disjoint from $A^sL_0$ over $A(t)$, since $\sqrt{1+t^2}\notin
L_0(i)$. \\[0.05in]
Let $L$ be the field composite of all
$L_n$'s, $n\in Q$ (if
$Q=\emptyset$,  we let $L=A(t)$). Then the extensions $A^s(t)$, $L_0$ and $L$ of $A(t)$
are all linearly disjoint over $A(t)$, and Galois, so that
$$\gal(A^sLL_0/A(t))\simeq G(A)\times \gal(L/A(t))\times
\gal(L_0/A(t))\simeq \hat\zee\times \gal(L_0/L)\simeq \hat\zee\times (\zee/2\zee)^\kappa.$$
Let $f_\alpha(X,Y)$, $\alpha<\kappa$, be an enumeration of
all absolutely irreducible polynomials of $L_0[X,Y]$ which are
separable in $Y$. We will construct by induction on $\alpha<\kappa$ a
chain $M_\alpha$ of algebraic extensions of $L_0$, which intersect $LA^s$
in $A(t)$ (and therefore are regular over $A$),  such that each
$M_{\alpha+1}$ is generated over $M_\alpha$ by a solution of
$f_\alpha(X,Y)=0$.  We let $M_0=L_0$,  and when $\alpha$
is a limit ordinal, we
let $M_\alpha=\bigcup_{\beta<\alpha}M_\beta$. Assume $M_\alpha$ already
constructed. \\[0.05in]
{\bf Claim}. $M_\alpha$ is Hilbertian. \\
Our assumption on $[M_{\beta+1}:M_\beta]$ being finite for $\beta<\alpha$ implies
that $M_\alpha=M'_\alpha L_0$, where $M'_\alpha$ is the union of $|\alpha|$ many finite algebraic
extensions of $A(t)$; hence the Galois closure $\tilde M_\alpha$ of
$M'_\alpha$ (over
$A(t)$) is the union of $|\alpha|$ many finite
Galois extensions of $A(t)$. Since $|\alpha|<\kappa$,  $L_0$ is not
contained in $\tilde M_\alpha$, and  $\gal(L_0/A(t)\simeq
(\zee/2\zee)^\kappa$ and Property \ref{hilb1}(2) give the result. \\[0.05in]
Construction of  $M_{\alpha+1}$.  We let $d$ be the degree of
$f_\alpha$  in $Y$, and
let $N_\alpha$ be the composite of all algebraic extensions of $L_0$ of
degree $\leq d$ and contained in $A^sL_0L$. Then $N_\alpha$ is a finite
(Galois)
extension of $L_0$. As $f_\alpha$ is absolutely
irreducible, it remains irreducible over $M_\alpha N_\alpha$. Because $M_\alpha$
is Hilbertian, there is some $a_\alpha\in M_\alpha$ such that
$f(a_\alpha,Y)$ is irreducible over $M_\alpha N_\alpha $ (by Property \ref{hilb1}(6)). We
let $M_{\alpha+1}$ be generated over $M_\alpha$ by a root $b_\alpha$ of
$f_\alpha(a_\alpha,b_\alpha)=0$. Since $[M_{\alpha+1}:M_\alpha]\leq d$,
$M_\alpha\cap A^sLL_0=L_0$,
and $M_{\alpha+1}$ is linearly disjoint from $N_\alpha M_\alpha$ over
$M_\alpha$, it follows that $M_{\alpha+1}\cap A^sL_0L=L_0$. \\[0.05in]
We let $M_\kappa=\bigcup_{\alpha<\kappa}M_\alpha$. By construction, every absolutely irreducible polynomial $f(X,Y)\in
L_0[X,Y]$ has a zero in $M_\kappa$, and therefore $M_\kappa$ is
PAC (by Fact \ref{rem1}(3)). \\
Recall that $L$ is linearly disjoint from $A^s$ over $A$, and that
$\gal(L/A(t))$ is procyclic.
Hence the topological generator $\si$ of $G(A)$ lifts to
an element $\si_2\in\gal(A^sL/A(t))$ whose restriction to $L$
topologically generates $\gal(L/A(t))$. As $M_\kappa$ is linearly
disjoint from $A^sL$ over $L_0$, this $\si_2$ lifts to an element
$\si'\in G(M_\kappa)$. Then $\langle\si'\rangle\simeq \hat\zee$ by Fact~\ref{rem1}(2),
$\fix(\si')$ is a regular extension of $A$ and is PAC. Hence the perfect
hull of $\fix(\si')$, $F_X$,  is pseudo-finite. \\[0.05in]
For the last assertion, note that $M_\kappa$ contains $L_0$, and that by
construction of $L_0$ and $L$, we have that $\cals(t,F_X)=\cals(t,M_\kappa)=X$
when char$(A)\neq 2$, and $\cals(t,F_X)\cap Z=X$ when
char$(A)=2$.

\begin{rem}\label{rem-prop1} Observe that if $A$ is infinite, and $F_X$ is as above, then
  the set of $Y$ such that $F_X\simeq _AF_Y$ has cardinality $\leq |A|$,
  since $|F_X|=|A|$. In particular, there are $2^{|A|}$ non-isomorphic
  models of $T(A)$ of the form $F_X$.

  \end{rem}

\begin{thm}\label{thm1} Let $T$ and $A$ be as above, with $A$ not
  pseudo-finite. Then $T(A)$ has no prime model.
\end{thm}

\prf Let us first do the very easy case when $A$ is finite. Then $T(A)$
is countable, and the existence of a prime model would imply that
isolated types are dense. If $|A|=q$, then  every model of $T(A)$ must
contain elements which are transcendental over $A$. In particular, by
\ref{psf1}(4), the formula $x^q\neq x$ contains no isolated type over
$A$. \\[0.05in]
Let us now assume that $A$ is infinite, and char$(A)\neq 2$. By
Proposition \ref{prop1}, a prime model $F$ of $T(A)$ has to (elementarily) $A$-embed in all
$F_X$'s, and therefore have  transcendence degree $1$ over $A$. Then
$F\prec F_X$ implies $F=F_X$. However, the set
$\cals(F)=\{\cals(u,F)\mid u\in F\setminus A\}$ has size $|A|$, hence
there is some subset $Y$ of $A$ which does not appear in
$\cals(F)$. I.e., $F$ cannot $A$-embed elementarily in that $F_Y$. So,
no prime model of $T(A)$ exists. A completely analogous discussion gives
the result in
 characteristic  $2$.

 \begin{rem} \label{rem-thm1}
   The proof of Theorem \ref{thm1} when $A$ is infinite only used item (2) of Proposition~\ref{prop1}. The
  interest of the first item is its formulation and relation to the
  following result of Jarden (see Theorems 18.5.6 and 18.6.1 in \cite{FJ}):
\begin{quote} {\em Let $K$ be a countable Hilbertian field. Then for almost
  all $\si$ in $G(K)$, the subfield of $K^{alg}$ fixed by $\si$ is
  pseudo-finite.}
\end{quote}
So, applying this to $K=A(t)$, we get that for almost all $\si$ in
$G(A(t))$, the subfield of $K^{alg}$ fixed by $\si$ is
  pseudo-finite. But the set of $\si$ with fixed subfield a regular
  extension of $A$ has measure $<1$ if  $G(A)\neq 1$, and for instance
  if $A=\ffi_p$, it has measure $0$. Item (1) of Proposition \ref{prop1}
  is therefore the correct generalisation: once fixed a lifting of a
  generator of $G(A)$, its tranlates by almost all elements of
  $G(A^s(t))$ fix a regular extension of $A$ which is
  pseudo-finite. Note that by Theorem 18.8.8 of \cite {FJ}, this result is false when $A$ is  uncountable.

\end{rem}

\section{Non-existence of prime saturated models}

\begin{defn} Let $T$ be a complete theory, $M$ a model of $T$.
\begin{enumerate}
\item The model $M$ is $\aleph_\varepsilon$-saturated if whenever $A\subset M$ is
  finite, then every strong type over $A$ is realised in $M$. When
  $T=T^{eq}$, equivalently, for any finite $A\subset M$, any type over
  $\acl(A)$ is realised in $M$.
\item Let $\kappa$ be an infinite cardinal or $\aleph_\varepsilon$. We say
  that $M$ is {\em $\kappa$-prime} if $M$ is
  $\kappa$-saturated, and elementarily embeds into every
  $\kappa$-saturated model of $T$.
\item Let $A\subset M$, and $T(A):={\rm Th}(M,a)_{a\in A}$, $\kappa$ as
  in (2). We say that
  $N$  is {\em $\kappa$-prime over $A$} if $N$ is a
  $\kappa$-saturated model of $T(A)$, and elementarily embeds into every
  $\kappa$-saturated model of $T(A)$.
\end{enumerate}
\end{defn}

\noindent
{\bf In the remainder of this section we assume GCH}.
\begin{rem}  It
  follows that if $\kappa$ is a regular cardinal larger than the
  cardinality of the language, and
  $A\subset M$ has cardinality $<\kappa$, then $T$ has $\kappa$-prime
  models over $A$, and furthermore they are all $A$-isomorphic: this
  follows easily observing that $T(A)$ has cardinality $<\kappa$, and
  the fact that our hypothesis on $\kappa$ guarantees that there are
  saturated models of $T(A)$ of cardinality $\kappa$. We will now show
  that when $T$ is the theory of a pseudo-finite field, this is
  essentially the only case when  $\kappa$-prime  models exist.
 The GCH hypothesis could be weakened to $2^{\aleph_0}=\aleph_1$ when
 $\trdeg(A)=\aleph_0$, and to $\lambda^{<\kappa}\leq \lambda^+$, where $\lambda=|A|$.
\end{rem}

%\newpage
\para\label{strategy}{\bf Setting and strategy}.  We let $\kappa$ be an
uncountable cardinal or $\aleph_\varepsilon$.  Let $A$ be a perfect
field of cardinality $\lambda\geq \kappa$, of infinite transcendence
degree with
absolute Galois group procyclic, and let $T(A)$ be the $\call(A)$-theory
whose models are the pseudo-finite fields which are regular extensions
of $A$. We also assume that $A$ is
not {a $\kappa$-saturated model of $T(A)$}, and fix a
transcendence basis $Z$ of $A$. Given a model $F$ of $T(A)$ and $t\in
F$, (almost) as before we define
$$\cals(t,F)=\begin{cases} \{a\in Z\mid \sqrt{t+a}\in F\} \hbox{ if }
  {\rm char}(A)\neq 2,\\
\{a\in Z\mid F\models \exists y\, y^2+y=at\} \hbox{ if }{\rm char}(A)=
2.
\end{cases}
$$
We will show that given a subset $X$ of
$Z$ with $|X|=\lambda$, there is a $\kappa$-saturated
model $F_X$ of $T(A)$ of cardinality {$\lambda^+$} and
with the following property: \\[0.05in]
{$(*)$\  For all $t\in F_X\setminus A$, there are some $b,c\in X$ such that
$b\in\cals(t,F_X)$ and $c\notin \cals(t,F_X)$.} \\[0.05in]
Note that this implies that  $\cals(t,F_X)\neq X$ and
$\cals(t,F_X)\neq Z\setminus X$ {for all $t\in F_X$}.

\medskip
Assume by way of contradiction that $F$ is a $\kappa$-prime model of
$T(A)$. Then it embeds elementarily in all fields $F_X$ constructed
above. Choose
$t\in F\setminus A$, and
let $Y=\cals(t,F)$. If $Y$ has size $\lambda$, then $F$ cannot
elementarily $A$-embed into $F_Y$. If $|Y|<\lambda$, then $Z\setminus Y$ has
size $\lambda$, and $F$ cannot elementarily $A$-embed into
$F_{Z\setminus Y}$ either. This shows that there is no $\kappa$-prime
model of $T(A)$.   \\[0.05in]
The proof of the theorem needs a technical lemma. Recall our setting:
$A$ is perfect of  cardinality $\lambda\geq \kappa$, with procyclic absolute Galois
group $G(A)$, and with transcendence basis $Z$.

\begin{lem} \label{tool} { Let $M$ and $C$ be perfect fields, both
    with procyclic absolute Galois groups, with
    $M/A$ regular, and with $M$ and $C$ linearly disjoint over $M\cap
    C$,  and regular
    over  $M\cap C$. We assume that $|M|=\lambda$, and
    that $|C|<\lambda$. We also fix a subset $X$ of $Z$ of cardinality
    $\lambda$. }
 %%  is relatively algebraically closed in $F^*$,  of cardinality
%%   $\lambda$.
  Then there is a regular extension $N$ of $M$ and of $C$, with
procyclic absolute Galois group, contained in $(MC)^{alg}$, and which
satisfies $(*)$ over $M$: if $a\in N\setminus M$, then there are $b,c\in
X$ such that $c\in \cals(a,N)$, but  $b\notin \cals(a,N)$.\end{lem}

\prf  We fix an enumeration $(b_\gamma)_{\gamma<\lambda}$, of
$(MC)^{alg}\setminus M$.  We will
build $N$ as the union of an increasing chain
$(MC_\gamma)_{\gamma<\lambda}$ of subfields of $(MC)^{alg}$, satisfying the following conditions:
\begin{itemize}[noitemsep]
\item[(i)] $C_{0}=C$.
\item[(ii)] If $\gamma$ is a limit ordinal, then
  $C_\gamma=\bigcup_{\delta<\gamma}C_\delta$.
\item[(iii)] If $\gamma<\lambda$, then
  $|C_\gamma|<\lambda$.
\item[(iv)] Each $C_{\gamma+1}$ is a  regular extension of $C_\gamma$, of
  finite transcendence degree over $C_\gamma$,  and
  $G(C_{\gamma})$ is procyclic.
\item[(v)] Each $MC_{\gamma}$ is a regular
  extension of $C_\gamma$ and of $M$.
  \item[(vi)] For each $\gamma$, $M$ and $C_\gamma$ are linearly disjoint
  over $M\cap C_\gamma$, and regular over $M  \cap C_\gamma$.
\item[(vii)]  \underline{Case (a)}: For each $\gamma$, either every $\si\in G(MC_\gamma)$ with
  $\fix(\si)$ regular over $M$ and over $C_\gamma$ moves $b_\gamma$, and
  in that case $b_\gamma\notin C_{\gamma+1}$. Or \\
\underline{Case (b)}: there is some $\si\in
  G(MC_\gamma)$, with fixed field regular over $M$ and over $C_\gamma$,
  and which contains $b_\gamma$; in that case $b_\gamma\in C_{\gamma+1}$,
and there are some elements $b,c\in X\cap C_{\gamma+1}$ such that $c\in
\cals(b_\gamma,C_{\gamma+1})$, and
$b\notin\cals(b_\gamma,C_{\gamma+1})$, and $b\in C_{\gamma+1}$.
\end{itemize}
It will then follow that
$N:=\bigcup_{\gamma<\lambda}MC_\gamma$ is a regular
extension of $M$ and of $C$, by (v). As it is regular over each
$C_\gamma$, and because every element of
$\bigcup_{\gamma<\lambda}MC_\gamma\setminus M$ occurs as a $b_\delta$,
it follows that $N$ satisfies $(*)$, by (vii) (and Remark \ref{rem2}). \\[0.05in]
{\bf Claim 1}. If $N$ is as constructed above, then $G(N)$ is
procyclic. \\[0.05in]
{\em Proof of Claim 1}.
Write $C_\lambda$ for $\bigcup_{\gamma<\lambda}C_\gamma$. Recall that
$M$ and $C_\lambda$ are linearly disjoint over their intersection, so
that $G(N)\simeq G(M)\times _{G(M\cap C_\lambda)}G(C_\lambda)$.
 If $G(N)$ is not
 procyclic, then there is a finite cyclic Galois extension $N_0$ of $M\cap
 C_\lambda$, and finite cyclic Galois extensions $N_1$ of $M$, and $N_2$ of
 $C_\lambda$, and positive integers $m$ dividing $n$, such that
 $[N_0:M\cap C_\lambda]=m$, $[N_1:M]=[N_2:C_\lambda]=n$, and
 $$\gal(N_1N_2/N)\simeq \gal(N_1/M)\times_{\gal(N_0/M\cap
   C_\lambda)}\gal(N_2/C_\lambda)\simeq
 \zee/n\zee\times_{\zee/m\zee}\times \zee/n\zee.$$ Thus $\gal(N_1N_2/N)$ has a cyclic subgroup $H$ of order $n$,
 projecting onto $\gal(N_1/M)$ and onto $\gal(N_2/C_\lambda)$, so that
 $\fix(H)$ is a proper Galois extension of $MC_\lambda$, which is regular over
 $M$ and over $C_\lambda$, contradicting item (vii) of the
 construction, since any generator of $\fix(H)$ over $MC_\lambda=N$
 appears as a $b_\delta$ for some $\delta<\lambda$. \qed

\medskip\noindent
We now start with the construction of the $C_\gamma$'s.
It is done by induction on $\gamma$, and if $\gamma$ is a limit ordinal, then we set
$C_\gamma=\bigcup_{\delta<\gamma}C_\delta$. Assume  that $C_\gamma$
has been constructed, we will now construct $C_{\gamma+1}$. If
$b_\gamma$ satisfies case (a) of (vii), then we let $C_{\gamma+1}=C_\gamma$. \\
Assume now that
we are in case (b) of (vii), and let $\si\in G(MC_\gamma)$ be such that $\fix(\si)$ is regular over $M$ and over
$C_\gamma$ and fixes $b_\gamma$.
Let  $Y_0\subset M$ be a finite set of algebraically
independent transcendentals (over $C_\gamma$),   such that
$b_\gamma\in C_\gamma(Y_0)^{alg}$, and let $b,c\in X$ be transcendental
and  algebraically independent over $C_\gamma(Y_0)$. We will construct
$C_{\gamma+1}$ as an algebraic extension of   $C_\gamma(Y_0,b_\gamma,b,c)$.
Let
$D=C_\gamma(Y_0)^{alg}\cap M=\bigl((M\cap
 C_\gamma)(Y_0)^{alg})\bigr)\cap M$,
$D_1=C_\gamma(Y_0,b,c)^{alg}\cap
M=\left((M\cap C_\gamma)(Y_0,b,c)\right)^{alg}\cap M$. \\[0.05in]
{\bf Claim 2}. $C_\gamma D(b_\gamma)$ and $M$ are linearly disjoint
over $D$, and $C_\gamma D(b_\gamma)$ is a regular extension of
$D$. The same hold for $C_\gamma D_1(b_\gamma)$ and $M$ over $D_1$.\\[0.05in]
{\it Proof of Claim 2}. Our hypothesis
on $b_\gamma$, and $Y_0\subset D\subset M$, imply the linear disjointness of
$C_\gamma D$ and $M$ over $D$, and  the second assertion. Moreover,
$M/D$ is regular, and $C_\gamma D(b_\gamma)$ is free from $M$ over
$D$. Hence, by Fact \ref{rem1}(8), $C_\gamma D(b_\gamma)$ and $M$ are
linearly disjoint over $D$. That $D\subset D_1\subset M$ gives the last
assertion. \qed

%% {\bf Claim 2}.
%% $C_\gamma D(b_\gamma)$ and $M$ are linearly disjoint over
%% $D$, and $C_\gamma D(b_\gamma)$ is a regular extension of
%% $D$.  \\[0.05in]
%% {\em Proof of Claim 2}. As $M$ and $C_\gamma$  are
%% linearly disjoint over $M\cap C_\gamma$, and $D\subset M$, it follows
%% that $C_\gamma D$ and $M$ are  linearly disjoint over $D$; moreover, by
%% Fact \ref{rem1}(7), as $C_\gamma(Y_0)^{alg}/C_\gamma(Y_0)$ is normal,
%% $Y_0\subset M$,
%% and $M$ is perfect, $C_\gamma(Y_0)^{alg}$ is linearly disjoint
%% from $M$ over $D$, and so is $C_\gamma D(b_\gamma)$, which proves the first
%% assertion. The second assertion follows from $MC_\gamma (b_\gamma)$ being regular
%% over $M$, and $M$ being regular over $D$. \qed

%\vspace{0.05in]
\medskip
\noindent
Let $E=C_\gamma D(b_\gamma, b,c,\sqrt{b_\gamma+b},\sqrt{b_\gamma+c})$ if
char$(A)\neq 2$, $E=C_\gamma D(b_\gamma,b, c, d_1,d_2)$, where
$d_1^2+d_1=bb_\gamma$ and $d_2^2+d_2=cb_\gamma$ if char$(A)=2$.  By
Corollary \ref{cor0}, as $C_\gamma D(b_\gamma)$ is a regular extension
of $D$, and is linearly disjoint from $M$ over $D$, it follows that $E$ is
regular over $C_\gamma D(b_\gamma)$ and over $D(b,c)$. By
Claim 2, this implies that $ME$ is regular over $M$ and over $C_\gamma
D(b_\gamma)$, and therefore also over $C_\gamma$.

$C_\gamma(Y_0)^{alg}$ are free and linearly disjoint over $D^{alg}$,
%$(C_\gamma(Y_0)^{alg}\cap M$,
we know that $E$ is regular over $D$ and over
$C_\gamma(Y_0, b_\gamma)$. Hence $ME$ is regular over $C_\gamma$ and $M$, and
$M$ and $D_1E$ are linearly disjoint over $D_1$. \\[0.05in]
Therefore, if $\si$ is a
topological generator
of $G(M)$, and $\si_\gamma$ is a topological generator of
$G(C_\gamma)$ which agrees with $\si$ on $(M\cap C_\gamma)^{alg}$, then we
may lift $\si_\gamma$ to an element $\si_{\gamma+1}$ of $G(C_\gamma
D (b_\gamma,b,c))$ which extends $\si$ on $D_1(b,c)^{alg}$,  which is the identity on $\sqrt{c+b_\gamma}$
and moves $\sqrt{b+b_\gamma}$ if char$(F)\neq 2$, and is the identity on
$d_2$ and moves $d_1$ if char$(F)=2$. We then let $C_{\gamma+1}$ be the
subfield of $C_\gamma D(b_\gamma,b,c)^{alg}$ fixed by $\si_{\gamma+1}$.
Then $M C_{\gamma+1}$ is a regular extension of $C_\gamma$ and of $M$. And by
construction, $\cals(b_\gamma,C_{\gamma+1})$ does not contain $b$ and
contains $b$. Moreover, by definition of $\si_{\gamma+1}$, any $M'$
containing $MC_{\gamma+1}$ and such that $b\in\cals(b_\gamma,M')$ will
contain $\sqrt{-1}$ if  char$(A)\neq 2$, and $\alpha$ if char$(A)=2$. I.e., $b\notin^*\cals(b_\gamma,MC_{\gamma+1})$. \qed

\begin{rem}\label{tool2} Let $A,C,M,N$ be as in Lemma \ref{tool}, and
  suppose that the $\lambda^+$-saturated model $F^*$ of $T(A)$ is a regular extension of $M$. Then there is an
  $M$-embedding $\varphi:N\to F^*$, with $F^*/\varphi(N)$
  regular.
  \end{rem}

\begin{proof} Since $N$ is a regular extension of $A$ with procyclic
  Galois group, we know that it embeds into a model of $T(A)$ which is
  regular over $N$. Consider $\tp(C/M)$: we know that $|M|=\lambda$,
  $|C|<\lambda$, so that by $\lambda^+$-saturation of $F^*$, $\tp(C/M)$
  is realized in $F^*$, say by $C'$; hence there is an $M$-isomorphism
  $\varphi:N\to (C'M)^{alg}\cap F^*$ which sends $C$ to $C'$ and with
  $F^*/\varphi(N)$ regular. \end{proof}

\begin{thm}\label{thm2} (GCH) Let $A$ and $T(A)$ be as in  \ref{setting}. Let
  $\kappa$ be an uncountable  cardinal, and suppose
  that $A$ is not a $\kappa$-saturated model of $T(A)$, and that
  $\kappa\leq |A|$. Then $T(A)$ has
  no $\kappa$-prime model. \\
If the transcendence degree of $A$ is $\geq \aleph_0$ and   $A$ is not $\aleph_\varepsilon$-saturated,
  then $T(A)$ has no $\aleph_\varepsilon$-prime model.
\end{thm}

\prf Let $\lambda=|A|$, $Z$ a transcendence basis of $A$, and let
$X\subset Z$ of size $\lambda$. We let $\kappa$ denote an uncountable
cardinal or $\aleph_\varepsilon$.
We fix a $\lambda^+$-saturated model $F^*$ of $T(A)$. We will work both
in $F^*$ and in some large algebraically closed field $\Omega$ containing  ${F^*}$.  We will build a
$\kappa$-saturated model $F_X$ of
$T(A)$ satisfying the  condition given in the strategy \ref{strategy}:
\begin{quote} $(*)$:
  For every $a\in F_X\setminus A$, $\cals(a,F_X)$ does not
  contain $X$ nor is it contained in $Z\setminus X$.
 \end{quote}
Then as explained above, the existence of these models will imply the
non-existence of $\kappa$-prime
models over $A$. While realising types within $F^*$ to obtain a
$\kappa$-saturated submodel of size $\lambda^+$ is easy,  condition
$(*)$  requires some work.\\
Recall that by Corollary~3.3 of \cite{H-psf}, if
$B$ is a relatively algebraically closed subfield of a pseudo-finite field $F$
such that $G(B)\simeq \hat\zee$, then the $\call(B)$-theory ${\rm Th}(F,a)_{a\in
  B}$ eliminates imaginaries. In that case, it follows that strong types
over $B$ are simply types over $B=B^{alg}\cap F$. We will
therefore start our construction by defining $F_0=A$ if $G(A)\simeq
\hat\zee$, and  if $G(A)\not\simeq \hat\zee$, then using Proposition \ref{prop1}(1) we  first choose some
 relatively
 algebraically closed
 subfield $F_{-1}$ of $F^*$ of transcendence degree $1$ over the prime
 field and such that
 the map $G(F^*)\to G(F_{-1})$ is an isomorphism. We may assume that
 $F_{-1}$ and $A$ are linearly disjoint over the relative algebraic
 closure of the prime field. We then apply Lemma
 \ref{tool} to the pair $(F_{-1}, A)$ to obtain a subfield
 $F'\subset\Omega$ of $(F_{-1}A)^{alg}$, satisfying the conclusion of
 Lemma \ref{tool}, and which we then move using Remark
 \ref{tool2} to a relatively algberaically closed subfield $F_0$ of $F^*$, regular
 over $A$, with $G(F_0)\simeq\hat\zee$, and satisfying condition ($*$)
 over $A$.

 \smallskip So, in both cases we have $G(F_0)\simeq
 \hat\zee$, and the same holds for all relatively algebraically closed subsets of $F^*$ containing
 $F_{-1}$. We will construct $F_X$ as a chain of $\lambda^+$ subfields
 of $F^*$. The reason for taking $\lambda^+$ instead of $\lambda$ is two-fold: First of all,
 $\lambda^{<\kappa}\leq \lambda^+=\lambda^\lambda$; and second, $\lambda^+$ is
 regular. \\[0.1in]
We use a diagonal argument, and build, by induction on
$\alpha<\lambda^+$, an increasing sequence $F_\alpha$ of subfields of
$F^*$, together with a collection of types
$(p_{\alpha,\beta})_{\beta<\lambda^+}$. To ease the writing, let us call
a subset $B$ of $F^*$ {\em small} if it is relatively
algebraically closed in $F^*$, has transcendence degree $<\kappa$ if
$\kappa\geq\aleph_1$, and finite if $\kappa=\aleph_\varepsilon$. We
choose the sequence $(F_\alpha)_{\alpha<\lambda^+}$ satisfying the
following conditions:
\begin{itemize}[noitemsep]
  \item[(a)] $F_0$ is as above.
\item[(b)] $F^*$ is a regular extension of $F_\alpha$, and $|F_\alpha|=\lambda$.
 \item[(c)] If $\alpha$ is a limit ordinal,
   $F_\alpha=\bigcup_{\beta<\alpha}F_\beta$.
   \item[(d)] $(p_{\alpha,\beta})_{\beta<\lambda^+}$ enumerates all (finitary)
     types over small
      subsets of $F_\alpha$.
   \item[(e)] $F_{\alpha+1}$ contains realisations of $p_{\delta,\beta}$ for
     all $\delta,\beta\leq \alpha$.
     \item[(f)] If $a\in F_{\alpha+1}\setminus F_\alpha$, then there are some
       $b,c\in X$ such that $\cals(a,F_{\alpha+1})$ contains $b$ but not
       $c$.
\end{itemize}
Items (a), (c) and (d) are atraightforward. Items (e), (f) and (b)
follow from Lemma \ref{tool} and its remark. Indeed, suppose $F_\alpha$ constructed; we
will build an increasing sequence of subfields $(M_\beta)_{\beta\leq\alpha+1}$  of $F^*$,
satisfying
\begin{itemize}[noitemsep]
  \item[(a')] $M_0=F_\alpha$.
\item[(b')] $F^*$ is a regular extension of $M_\beta$, and $|M_\beta|=\lambda$.
 \item[(c')] If $\beta$ is a limit ordinal,
   $M_\beta=\bigcup_{\gamma<\beta}M_\gamma$.
    \item[(e')] If $\beta\leq\alpha$, then $M_{\beta+1}$ contains realisations of
      $p_{\beta,\alpha}$
     and  of $p_{\alpha,\beta}$.
     \item[(f')] If $a\in M_{\beta+1}\setminus M_\beta$, then there are some
       $b,c\in X$ such that $\cals(a,M_{\beta+1})$ contains $b$ but not
       $c$.
\end{itemize}
Item (a') and (c') are straightforward. Assume $M_\beta$ given, we
will construct $M_{\beta+1}$ as follows: let $E$ be a
small subset of $M_\beta$ containing the bases of $p_{\beta,\alpha}$ and
$p_{\alpha,\beta}$, choose $(a_1,a_2)$  realising
$(p_{\beta,\alpha},p_{\alpha,\beta})$ in $F^*$, transcendental and  algebraically
independent over $M_\beta$. Now apply Lemma~\ref{tool} and its
Remark~\ref{tool2} to
$(C,M)=(E(a_1,a_2)^{alg}\cap F^*, M_\beta)$ to obtain first $N$, then $\varphi(N)=M_{\beta+1}$
satisfying (b'), (e') and (f'). {Note that in applying Remark
  \ref{tool2}, we will take for $\varphi$ an $M_\beta$-embedding of $N$ into $F^*$.} \\[0.05in]
We then let
$F_{\alpha+1}=\bigcup_{\beta\leq \alpha+1}M_{\beta}$. Then item (b) and
(f) hold (because at  stage $\alpha+1$, we are only realising $\lambda$ many types);
for (e), let $\beta,\delta\leq \alpha$. If  $\beta,\delta<\alpha$,
then $p_{\beta,\delta}$ is realised in $F_\alpha$. If  $\beta=\alpha$, then
$p_{\beta,\delta}$ is realised in  $M_{\delta+1}$, and if
$\delta=\alpha$, in $M_{\beta+1}$, which shows (e).\\[0.1in]
Define $F_X=\bigcup_{\alpha<\lambda^+}F_\alpha$. By construction, $F_X$ is
$\kappa$-saturated: if $q$ is a type over the small set
$B\subset F_X$,
then $B\subset F_\alpha$ for some $\alpha<\lambda^+$, and therefore $q$
appears as a $p_{\alpha,\beta}$ for some $\beta<\lambda^+$. Then $q$ is
realised in $F_\gamma$, where $\gamma=\sup\{\alpha,\beta\}+1$. In
particular, $F_X$ is PAC, perfect, with absolute Galois group isomorphic to
$\hat\zee$, and therefore pseudo-finite. \\[0.1in]
Furthermore, if $a\in F_X\setminus A$, then $a\in F_\alpha$ for some
$\alpha<\lambda^+$, and therefore
$\cals(a,F_\alpha)=\cals(a,F_X)$ neither contains $X$ nor is contained in
$X\setminus F$. This finishes the proof of the Theorem. \qed

%%  $$\begin{diagram}
%%  \node{M_\beta}\arrow{s,-}\\
%%  \node{D_\gamma(b,c)}\arrow{s,-}\\
%%  \node{D_\gamma}\arrow{e,-}\arrow{s,-}\node{C_\gamma
%%    D_\gamma}\arrow{e,-}\arrow{s,-}\node{C_\gamma D_\gamma(b_\gamma)}\\
%%  \node{(C_\gamma\cap
%%    M_\beta)(X_0)}\arrow{e,-}\arrow{s,-}\node{C_\gamma(X_0)}\arrow{s,-}\\
%%  \node{C_\gamma\cap M_\beta}\arrow{e,-}\node{C_\gamma}
%%  \end{diagram}$$

%%  \node[5]{C_\gamma D_\gamma(b_\beta)}\arrow{ssw,-}\\
%%  \node[2]{D_\gamma(b,c)}\arrow{sse,-}\node[4]{C_\gamma
%%    D_\gamma}\arrow{ssw,-}\arrow{sse,-}\node[6]{C_\gamma D_\gamma(b_\beta)}\arrow{sse,-}\\
%%  \node[3]{D_\gamma}\arrow{sse,-}\node[5]{C_\gamma (X_0)}\arrow{sse,-}\arrow{ssw,-\\
%%  \node[4]{(C_\gamma\cap
%%    M_\beta)(X_0)}\arrow{sse,-}\node[6]{C_\gamma}\arrow{ssw,-}\\
%%  \node[5]{C_\gamma\cap M_\beta}

\noindent

\para{\bf Concluding remarks}. When $A$  has finite transcendence degree, then the above
  proof breaks down. It might be possible to fix it by choosing a
  suitable subset $Z$ of $A$. \\
In the absence of GCH, saturated models will
in general not exist. The construction given above works when
$\lambda^{<\kappa}\leq \lambda^+$ when $\lambda=|A|$ is uncountable, and
under CH when $\kappa=\aleph_\epsilon$ and $|A|$ has transcendence
degree $\aleph_0$, and  the same argument shows that
$T(A)$ has no $\kappa$-prime model.\\
I believe that Theorems \ref{thm1} and \ref{thm2} generalise to the case of arbitrary
bounded PAC
fields without assuming perfection (Recall that a field is bounded if it has finitely many Galois extensions
of degree $n$ for every $n>1$; Claim 1 of Lemma \ref{tool} will need a
different proof).   %% What will require some work is when
%% $G(F)\not\simeq G(A)$: one needs to generalise appropriately the proof
%% of Proposition \ref{prop1}(2) and find a  Galois extension $L$ of
%% $A(t)$, regular over $A$ and  such that $\gal(LA^s(t)/A(t))\simeq
%% G(F)$.
%% Proposition \ref{prop1}(1) generalises easily to $e$-free perfect
%% PAC fields.

%\end{rems}

\end{document}